\documentstyle[amssymb,11pt,epsf]{amsart}
\newtheorem{theorem}{Theorem}[section]
\newtheorem{proposition}[theorem]{Proposition}
\newtheorem{lemma}[theorem]{Lemma}

\newtheorem{corollary}[theorem]{Corollary}

\newenvironment{proof}{\trivlist\item[\hskip \labelsep {\sc
  Proof:}\enskip]}{\unskip\nobreak\hskip 2em plus 1fil\nobreak%
\hfill\fbox{\rule{0ex}{.5ex}\hspace{.5ex}\rule{0ex}{.5ex}}\endtrivlist}

\newcommand{\bqn}{\begin{equation}}
\newcommand{\eqn}{\end{equation}}
\newcommand{\R}{{\Bbb R}}
\newcommand{\C}{{\Bbb C}}
\newcommand{\Z}{{\Bbb Z}}
\newcommand{\F}{{\Bbb F}}

\renewcommand{\S}{{\Bbb S}}
\renewcommand{\sl}{\mbox{sl}}
\newcommand{\bd}{\partial}
\renewcommand{\H}{{\Bbb H}}

\newcommand{\cop}{{\mathop{cop}}}
\newcommand{\op}{{\mathop{op}}}

\newcommand{\free}{{\mbox{\scriptsize free}}}
\newcommand{\opname}[1]{\mathop{\fam0#1}}

\newcommand{\Ie}{{\it I.e.}}
\newcommand{\ie}{{\it i.e.}}
\newcommand{\Eg}{{\it E.g.}}
\newcommand{\eg}{{\it e.g.}}
\newcommand{\Ad}{\mbox{Ad}}
\newcommand{\SU}{\opname{SU}}
\newcommand{\SO}{\opname{SO}}
\newcommand{\Tr}{\mbox{Tr}}
\newcommand{\Hom}{\mbox{Hom}}

\newcommand{\End}{\mbox{End}}

\newcommand{\im}{\mbox{im\ }}

\newcommand{\tensor}{\otimes}

\newcommand{\longto}{{\longrightarrow}}
\newcommand{\manyto}{{{}^\longto_\longto}}
\newcommand{\action}{\to}

\newcounter{fignum}
\newcounter{arrownum}
\newcommand{\arepsf}[1]{\epsfbox{#1.arrow}}
\newcommand{\arrowpic}{$$\addtocounter{arrownum}{1}
    \arepsf{\thearrownum}$$}
\newcommand{\figepsf}[1]{\epsfbox{#1.fig}}
\newcommand{\fig}{$$\addtocounter{fignum}{1}
    \figepsf{\thefignum}$$}

\begin{document}

\title{Non-involutory Hopf algebras and 3-manifold invariants}
\author{Greg Kuperberg}
\address{Department of Mathematics \\ Yale University\\
        New Haven, CT 06520}
\email{greg@@math.yale.edu}
\thanks{The author was supported by an NSF Postdoctoral Fellowship,
        grant \#DMS-9107908.}
\date{July 11, 1995}

\begin{abstract}
We present a definition of an invariant $\#(M,H)$, defined for
every finite-dimensional Hopf algebra (or Hopf superalgebra or Hopf
object) $H$ and for every closed, framed 3-manifold $M$.  When $H$
is a quantized universal enveloping algebra, $\#(M,H)$
is closely related to well-known quantum link invariants such as
the HOMFLY polynomial, but it is not a topological quantum field
theory.
\end{abstract}

\maketitle

This paper presents a definition of an invariant $\#(M,H)$ which depends on a
framed, closed 3-manifold $M$ and a finite-dimensional Hopf algebra $H$, and
whose value lies in the ground field of $H$.  The Hopf algebra $H$ need not
be quasitriangular, triangular, ribbon, modular, a quantum deformation,
involutory, or semisimple, nor does it need to have any other decoration or
structural property.  It can be any finite-dimensional example of the object
defined by Sweedler \cite{Sweedler} and Drinfel'd \cite{Drinfeld} or, more
generally, a finite-dimensional Hopf super-algebra or a Hopf object in any
category which sufficiently resembles the category of finite-dimensional
vector spaces.  In a previous paper \cite{Kuperberg:hopf}, the author defined
$\#(M,H)$ for involutory Hopf algebras (Hopf algebras in which the square of
the antipode is the identity) and for closed and unoriented but unframed
3-manifolds.  An important intermediate class of finite-dimensional Hopf
algebras is the class of balanced Hopf algebras, for which the 3-manifold $M$
need only be oriented and combed rather than framed.  Recall that a framing
of a 3-manifold is a homotopy class of linearly independent triples of
tangent vectors fields.  A combing is defined as the homotopy class of a
single non-vanishing tangent vector field.

In a subsequent paper \cite{Kuperberg:copout}, we will define the related invariant
$\#(M,L,H)$, where $M$ is a framed, closed 3-manifold, $H$ is a Hopf algebra,
and $L$ is a framed link in $M$.  More generally, the invariant $\#(M,G,H)$
can be defined, where $G$ is a framed graph in $M$.  When $M = \S^3$, these
invariants coincide with the Reshetikhin-Turaev invariants of links and
ribbons graphs derived from $D(H)$, the quantum double of $H$.  In
particular, if $q$ is a root of unity and $\frak g$ is a simple Lie algebra,
the Reshetikhin-Turaev invariants for the finite-dimensional quantum groups
$u_q(\frak g)$ yield root-of-unity values of the familiar quantum link
invariants, such as the Jones polynomial, the HOMFLY polynomial, the Kauffman
polynomial, and the quantum $G_2$ link invariant.  The Hopf algebra
$u_q(\frak g)$ is almost the quantum double of $u_q(\frak g^+)$, a truncated
quantum deformation of (the enveloping algebra) $U(\frak g^+)$, where $\frak
g^+$ is a Borel subalgebra of $\frak g$.  Therefore $H = u_q(\frak g^+)$ is
an important special case of the invariant $\#(M,H)$ that we define here.

Some other important special cases of $\#(M,H)$ are the following: $\#(\S^3,H)
= 1$ by normalization, while $\#(\S^2 \times \S^1,H) = \Tr(S^2)$ is $\dim H$
when $H$ is involutory and 0 when $H$ is non-involutory, and $\#(\R P^3,H) =
\Tr(S)$.  (Here $S$ is the antipode map of $H$ and $\Tr$ is its trace as a
linear endomorphism of $H$. To distinguish the $n$th power of the antipode
from the $n$-sphere, we write the former as $S^n$ and the latter as $\S^n$.)
Moreover, $\#(M,H)$ is multiplicative under connected sums and under tensor
products of Hopf algebras. If $H = \C[G]$ is the Hopf algebra of a group,
$\#(M,H)$ is the number of homomorphisms from the fundamental group of $M$ to
$G$, which can be written as $|\Hom(\pi_1(M),G)| = |[M:B_G]| = |H^1(M,G)|$.
If $H$ is an exterior algebra with one generator, which is a two-dimensional
graded Hopf algebra, then $\#(M,H) = |H_1(M,\Z)|$ when the right side is
finite, and $\#(M,H) = 0$ otherwise.

Given the relation between $\#(M,L,u_q(\frak g^+))$ and quantum link invariants,
one might suspect, as the author once did, that $\#(M,u_q(\frak g^+))$ includes
or is equivalent to the Jones-Witten 3-manifold invariants, defined explicitly
by Reshetikhin and Turaev \cite{RT:manifolds} for the group $\SU(2)$ and by
Turaev and Wenzl \cite{Turaev-Wenzl} for the group $\SU(n)$. However, they cannot
be equal or equivalent up to normalization, because $\#(\S^2 \times \S^1, H) =
0$ when $q \ne 1$ (and in general when $H$ fails to be semisimple and
co-semisimple), which violates the axioms of a topological quantum field theory
(TQFT), and the relevant Hopf algebras are non-involutory.  Rather,
$\#(M,H)$ should agree with another invariant recently found by Kauffman and
Radford \cite{Kauffman-Radford} and considered by Reshetikhin
\cite{Reshetikhin} which has a Dehn surgery definition similar to the definition
of the Reshetikhin-Turaev 3-manifold invariant.

If $H$ is involutory and $\dim H \ne 0$, then $\#(M,H)$ is a TQFT. In this case
$H$ must be a semisimple algebra \cite{Larson-Radford}, and
$\#(M,H)$ can equally well be defined in Turaev-Viro style using the
representation category of $H$ \cite{Barrett-Westbury}, or in
Reshetikhin-Turaev-Wenzl style using the representation theory of $D(H)$.  In
the semisimple case, the Kauffman-Radford invariant also equals the
corresponding Reshetikhin-Turaev-Wenzl-style invariant.

\section{The basic idea; involutarity}
\label{s:invol}

Let $G$ be a finite group, not necessarily commutative, and let $C$ be a
connected, finite simplicial complex with oriented edges and faces and with a
distinguished vertex $p$ chosen as a base point.   Recall the definition of
$H^1(C,G)$, the non-commutative first cohomology of $C$ with coefficients in
$G$:  A 1-cochain is a function from the edges of $C$ to $G$. Given a cochain
and a face, let $g_1$, $g_2$, and $g_3$ be the three group elements assigned to
the edges of a face in the order given by the orientation of the face, and let
$\sigma_1$, $\sigma_2$ and $\sigma_3$ each be $1$ or $-1$ if the orientation of
the corresponding edge agrees or disagree with the orientation of the face:
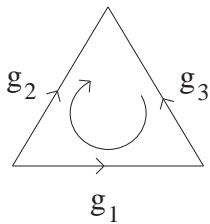
\fig
The cochain is a co-cycle if $$g_1^{\sigma_1} g_2^{\sigma_2} g_3^{\sigma_3} =
1$$ for every face; here $\sigma_1 = \sigma_3 = -1$ and $\sigma_2 = 1$.  A
0-cochain is a function from the vertices of $C$ to $G$, and a 0-cochain is
reduced if it is 1 on the base point $p$. The 0-cochains form a group under
pointwise multiplication with the reduced 0-cochains as a subgroup, and the
coboundary operation can be understood as a group action of 0-cochains
actings on 1-cochains:  Given a 0-cochain $c$ and a 1-cochain $d$, if $d$ has
the value $g$ at an edge $e$ and $c$ has the values $h_1$ and $h_2$ at the
head and tail of $e$, then the 1-cochain $cd$ has the value $h_1 gh_2^{-1}$
at $e$.  The cohomology set $H^1(C,G)$ is defined as the set of orbits of the
0-cochains acting on 1-cochains, while the reduced cohomology set
$\widetilde{H}^1(C,G)$ is defined as the set of orbits of the reduced
0-cochains. The reduced cohomology is the set of group homomorphisms from
$\pi_1(C)$ to $G$, while the full cohomology is the set of conjugacy classes
of these group homomorphisms.  Since reduced 0-cochains act freely on
1-cochains, the number of 1-co-cycles can be computed as
$|G|^{v-1}|\Hom(\pi_1(C),G)|$, where $v$ is the number of vertices of $C$.  

The above analysis fully extends from simplicial complexes to polygonal
complexes.  A {\em polygonal complex} is a CW complex whose 2-skeleton
consists of a collection of polygons whose edges are identified to each other
by homeomorphisms and whose vertices may be identified through further
equivalences.  A co-cycle is still defined as a $G$-valued function on edges such
that, up to orientations, the product of the values of the edges around every
face is the identity, and the conclusion that the total number of 1-co-cycles is
$|G|^{v-1}|\Hom(\pi_1(C),G)|$ still holds.  Removing the powers of $|G|$ by
normalization, define $\#(C,G)$ to be the remaining factor
$|\Hom(\pi_1(C),G)|$.  Clearly, $\#(C,G)$ is a topological invariant.

The first goal is to redefine $\#(C,G)$ in terms of the Hopf algebra $\F[G]$,
where $\F$ is a field of characteristic 0.  Recall that the Hopf algebra
structure on $\F[G]$ consists of an algebra structure which is the linear
extension of group multiplication, and an algebra structure on the dual vector
space $\F(G)$ which is pointwise multiplication of $\F$-valued functions on $G$.
Recall also that the antipode map $S:\F[G] \to \F[G]$ is the linear extension of
the group inverse.  Let the elements of $G$ be both a standard basis and
standard dual basis for $\F[G]$. Let $M_{a_1 a_2 \ldots a_n}$ be the tensor
corresponding to the multilinear form $\Tr(A_1 A_2 \ldots A_n)$, where the trace
is the trace of left (or right) multiplication, and let $\Delta^{a_1 a_2 \ldots
a_n}$ be the same tensor defined by the dual algebra structure.  Taking the
indices as group elements, observe that $M_{a_1 a_2 \ldots a_n}$ is $|G|$ if
$a_1 a_2 \ldots a_n = 1$ and is zero otherwise, while $\Delta^{a_1 a_2 \ldots
a_n} = 1$ if all $a_i$'s are equal and is zero otherwise.  The fact that the
matrix of an $M$ tensor expresses the co-cycle rule of non-commutative cohomology
can be exploited as follows:  Assign to each edge a $\Delta$ tensor with an
index for each incident face, and assign to each face an $M$ tensor with an
index for each edge.  Assign to each pair consisting of a face and an edge it
contains either the antipode map or the identity map from $\F[G]$ to $\F[G]$,
depending on whether the orientations disagree or agree.  Take the tensor
product of all tensors and contract according to incidence.  The resulting
element of $\F$ is evidently the sum of $|G|^f$ over all 1-co-cycles,
and therefore equals $|G|^{v+f-1}\#(C,G)$ if $C$ has $v$ vertices
and $f$ faces.

If a polygonal CW complex $C$ is a cell decomposition of an oriented 3-manifold
$M$, an interesting symmetry appears.  Consider the handle decomposition of $M$
corresponding to $C$, i.e., thicken each face and replace each edge by a prism
with a long, thin side for each incident face prism:
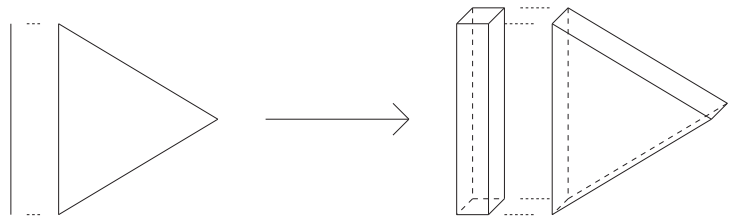
\fig
Call the long thin sides of either kind of prism {\em rectangles}.  The
orientation of $M$ is a right-hand rule, which associates a {\em
co-orientation} or a cyclic ordering of the rectangles of the prism of an
oriented edge.  Working from the handle decomposition, the tensorial
definition of $\#(C,G)$ associates a $\Delta$ tensor to each edge prism, an
$M$ tensor to each face prism, and either $I$ or $S$ to each rectangle.

If the indices of the $\Delta$ tensors are cyclically ordered according
to the co-orientations, the tensorial expression is evidently invariant
under an exchange of $\Delta$ and $M$ by Poincar\'e duality. As stated,
it is a well-defined formula for any finite-dimensional Hopf algebra
$H$.  It is a familiar topological invariant for a large class of
co-commutative Hopf algebras, namely group algebras, and it is the same
invariant for a large class of commutative Hopf algebras, namely dual
group algebras.  One might conjecture, as the author once did, that it
is always a topological invariant, but this is only true when the
antipode $S$ is an involution, \ie, when $H$ is involutory.  (See
reference~\cite{Kuperberg:hopf} for a proof of  invariance.  The
involutory invariant was also independently discovered by Chung,
Fukuma, and Shapere \cite{Shapere}.)

A valid generalization to non-involutory Hopf algebras is the subject of the
rest of the paper.  The definition of $\#(M,H)$ for $H$ involutory mandates
that a rectangle be replaced with $S^n$ (the $n$th power of the antipode),
where $n \in \Z/2$ is 0 if the link of the edge prism and the rim of the face
prism intersect positively and 1 if they intersect negatively on the surface
formed by all rectangles, which is an example of a Heegaard surface.  This
rule is not natural (and does not lead to an invariant) when the exponent $n$
is an integer rather than an element of $\Z/2$; we must choose between all
even $n$ or all odd $n$, depending on the sign of the orientation.  The
solution is to consider a logarthmic lifting of the unit tangent bundle of
the Heegaard surface to an affine line bundle, which allows angles to take
values in the real numbers rather than the circle and intersections to take
values in the integers rather than in $\Z/2$.  Unless the Heegaard surface is
a torus, such a lifting must be singular, but its singularities can solve
other problems that arise elsewhere in the definition.  A singular combing of
the Heegaard surface which extends to a combing of the 3-manifold can provide
just such a logarithmic structure.

\section{Topology}

Throughout the paper, a manifold is a compact, oriented, triangulated,
differentiable manifold.  In addition, a manifold is assumed to have no boundary
unless otherwise stated.  PL and smooth constructions will be freely combined, a
liberty which rarely leads to  confusion in three and fewer dimensions.

\subsection{Combings and framings}

Let $M$ be a 3-manifold.  A {\em combing} of $M$ is a non-vanishing
tangent vector field, and a {\em framing} of $M$ consists of three
linearly independent vector fields whose orientation agrees with that
of the manifold. For convenience, view $M$ as a Riemannian manifold, a
combing as a section of the unit tangent bundle, and a framing as three
orthogonal sections.  Since the orientation of $M$ induces a
cross-product operation on tangent vectors, it suffices to describe a
framing as a pair of orthogonal sections, the third section being given
by the cross product.

The goal of this subsection is to classify the set of combings and the set of
framings of a 3-manifold up to homotopy. It is important to describe these
two sets explicitly, since the invariant $\#(M,H)$ in general depends on a
combing or framing of $M$.  In addition, the proofs of Lemmas~\ref{lcrealize}
and \ref{lcequiv} depend on the classification.

It is  a classical result that the tangent bundle $TM$ of a 3-manifold $M$ is
trivial.  Using a trivialization of $M$, which is itself a framing, a combing
becomes a map from $M$ to $\S^2$ and a framing becomes a map from $M$ to
$\SO(3)$. Thus, the set of combings is bijective with the homotopy set
$[M,\S^2]$ and the set of framings is bijective with $[M,\SO(3)]$.  However,
for many $M$, since no one trivialization of $M$ is completely natural,
combings and framings are not canonically bijective with $[M,\S^2]$
and $[M,\SO(3)]$.  Rather, it is important to understand the natural
group action of $[M,\SO(3)]$ on $[M,\S^2]$ and $[M,\SO(3)]$.

\begin{proposition} The homotopy class of a map from $M$ to $\S^2$ is
described by a {\em characteristic class} $c \in H^2(M,\Z)$, defined as the
pullback of $[\S^2]$, and a {\em Hopf degree} $d$ in an affine space of
$\Z/n$, where $n$ is defined as follows:  Let $H^2_\free(M,\Z)$ be the
quotient of $H^2(M,\Z)$ by its torsion, and let $c_\free$ be the image of $c$
under this quotient.  Then $n$ is the {\em maximal divisor} of $c_\free$, the
largest integer such that $c_\free/n$ exists in the lattice $H^2_\free(M,\Z)$.
\end{proposition}

\begin{proof}(Sketch)
In general, if $X$ and $Y$ are CW complexes, the homotopy set $[X,Y]$ is
approximated by the direct product of the cohomology groups of $X$ with
coefficients in the homotopy groups of $Y$.  More precisely, if $X_n$ is the
$n$-skeleton of $X$, then restriction from $X_n$ to $X_{n-1}$ yields a map
$j_n:[X_n,Y] \to [X_{n-1},Y]$, and of course $[X,Y]$ is the inverse limit. 
Moreover, there is a group action of $H^n(X,\pi_n(Y))$ on $\im j_{n+1}$, and
each orbit of this group action is $j^{-1}(f) \cap \im j_{n+1}$ for some $f
\in [X_{n-1},Y]$.  In simple cases, each $j_n$ is a surjective group
homomorphism with a complemented kernel and the action of $H^n(X,\pi_n(Y))$ is
free, in which case $[X,Y] \cong H^*(X,\pi_*(Y))$.  But the homotopy set
$[M,\S^2]$ exhibits a complication that arises because $\S^2$ has non-trivial
homotopy groups of adjacent degree.  Although $\pi_2(\S^2) \cong \pi_3(\S^2)
\cong \Z$, $[M,\S^2]$ resembles but does not necessarily equal $H^2(M,\Z)
\oplus H^3(M,\Z)$.  To obtain the exact answer, consider the fibration
sequence
\begin{equation}
\S^1 \to \S^3 \to \S^2 \to \C P^\infty \to \H P^\infty, \label{e:hopffib}
\end{equation}
where the first three terms are the Hopf fibration.  The terms $\C P^\infty$
and $\H P^\infty$ are classifying spaces of $\S^1$ and $\S^3$ viewed as groups,
and they extend the Hopf fibration because their loop spaces are homotopy
equivalent to $\S^1$ and $\S^3$. Sequence~(\ref{e:hopffib}) would lead to the
exact sequence
$$[M,\S^1] \to [M,\S^3] \to [M,\S^2] \to [M,\C P^\infty] \to [M,\H P^\infty]$$
if all of the terms in the sequence were groups.  Although three of the terms
can be recognized as the cohomology groups of the 3-manifold $M$ and the last
term is trivial, the homotopy set $[M,\S^2]$ is not a group.  Instead,
corresponding exact sequence is
$$H^1(M,\Z) \manyto H^3(M,\Z) \action [M,\S^2] \to H^2(M,\Z) \to 0,$$
where the second map is a group action and the first map varies according to a
chosen orbit of the group action.  Recall that the Hopf degree establishes the
isomorphism $\pi_3(\S^2) \cong \Z$; a calculation in obstruction theory using
the Hopf degree establishes that the first arrow is:
$$H^1(M,\Z) \stackrel{\cup c}{\longto} H^3(M,\Z) \action [M,\S^2]
\to H^2(M,\Z) \ni c.$$
In other words, if an orbit of $[M,\S^2]$ is given by a cohomology class $c \in
H^2(M,\Z)$, the stabilizer in $H^3(M,\Z)$ of this orbit is the image of the cup
product with $c$.  Recall that $H^1(M,\Z)$  is free, and that the cup product
$$\cup:H^1(M,\Z) \times H^2(M,\Z) \to H^3(M,\Z)$$
annihilates the torsion of $H^2(M,\Z)$ and is a unimodular bilinear
form, or perfect pairing, on the free parts of the cohomology groups.
It follows that the image of $d \mapsto d \cup c$ is generated by
$n[M]$, where $n$ is the maximal divisor of $c_\free$.
\end{proof}

\begin{proposition} The homotopy set $[M,\SO(3)]$ is a central extension of
$H^1(M,\Z/2)$ by $H^3(M,\Z)$.  The group action of $f \in [M,\SO(3)]$ on
$[M,\S^2]$ is partly characterized as follows: The group $H^1(M,\Z/2) =
[M,\SO(3)]/H^3(M,\Z)$ has a quotient $H^2(M,\Z) * \Z/2$ by the universal
coefficient theorem; let $f'$ be the image of $f$ in $H^2(M,\Z) * \Z/2$.
quotient.  Then $f$ acts on the characteristic class $c$ of an element of
$[M,\S^2]$ by $c \mapsto c+f'$, if $H^2(M,\Z) * \Z/2$ is understood as a
subgroup of $H^2(M,\Z)$.  If $f \in H^3(M,\Z)$, then $f$ acts on the degree
$d$ of an element of $[M,\S^2]$ by $d \mapsto d+f$. \label{pso3action}
\end{proposition}
\begin{proof}(Sketch)
The homotopy set $[M,\SO(3)]$ is much easier to understand, since the
non-trivial homotopy groups $\pi_1(\SO(3)) \cong \Z/2$ and $\pi_3(\SO(3))
\cong \Z$ are not adjacent.  In general, $[M,\SO(3)]$ is a central extension
of $H^1(M,\Z/2)$ by $H^3(M,\Z)$ as a group, where the group law is defined
using the Lie group structure of $\SO(3)$.  In terms of the geometry on $M$,
one framing $f_1$ converts a second framing $f_2$ into an element of
$[M,\SO(3)]$, and after this conversion the $H^1(M,\Z/2)$ class of $f_2$
describes its spin structure, while the $H^3(M,\Z)$ class gives its degree.
Since these classes are defined relative to $f_1$, the set of homotopy
classes of framings of $M$ is only an affine space of the group $[M,\SO(3)]$.
However, the action of $[M,\SO(3)]$ on $[M,\S^2]$ is again more complicated. 
There is a fibration sequence analogous to sequence~(\ref{e:hopffib}):
\begin{equation}
\S^1 \to \SO(3) \to \S^2 \to \C P^\infty \to B_{\SO(3)}. \label{e:so3fib}
\end{equation}
Here is an explicit description of the third map:  $\SO(3)$ is a principal bundle
over $\S^2$ with fiber $\S^1$, and since $\C P^\infty$ is the classifying space of
the group $\S^1$, the bundle over $\S^2$ induces a map $\S^2 \to \C P^\infty$ which
is the map in the sequence.  In sequence~(\ref{e:hopffib}), the
third map has degree 1, because the Hopf fibration has Chern class 1.  By
contrast, since $\SO(3)$ has Chern class 2, the third map of
sequence~(\ref{e:so3fib}) has degree 2.  Also, recall the familiar fact that the
inclusion $\S^1 \to \SO(3)$ is a mod 2 reduction $\Z \to \Z/2$ of $\pi_1$.

Sequence~(\ref{e:so3fib}) yields an exact sequence of homotopy sets
of maps from $M$, which simplifies to another mixed exact sequence of sets and
groups:
$$H^1(M,\Z) \to H^1(M,\Z/2) \oplus H^3(M,\Z) \to [M,\S^2] \to
H^2(M,\Z) \to H^2(M,\Z/2).$$
The action of the $H^3(M,\Z)$ subgroup of $[M,\SO(3)]$ on $[M,\S^2]$ can
be understood by recalling that the covering map and group homomorphism
$\SU(2) \to \SO(3)$ is an isomorphism of $\pi_3$, and, taking $\S^3
\cong \SU(2)$, the group action of $[M,\SU(2)]$ on $[M,\S^2]$ induced
by this homomorphism is the same as the one in
sequence~(\ref{e:hopffib}).  So this part of the action only affects
the Hopf degree $d$.

Restricting attention to the $H^2(M,\Z)$ quotient of $[M,\S^2]$ given by the
characteristic class, the exact sequence further simplifies to
\begin{equation} H^1(M, \Z) \stackrel{\bmod 2}{\longto} H^1(M,\Z/2)
\stackrel{\phi}{\longto} H^2(M,\Z) \stackrel{\times 2}{\longto} H^2(M,\Z)
\stackrel{\bmod 2}{\longto} H^2(M,\Z/2), \label{e:twodim}\end{equation}
where the map $\phi$ is the action of $[M,\SO(3)]$ on the characteristic class
of $[M,\S^2]$.  Sequence~(\ref{e:twodim}) demonstrates that while the
characteristic class $c$ of a map to $\S^2$ is affected by the action of
$\SO(3)$, $2c$ is not. In terms of the geometry on $M$, the class $c$ of a
combing depends on the choice of framing, but $2c$ does not.   Since $2c$ is
the pullback of the Euler class $2[\S^2]$ on $\S^2$, it can be recognized as
the Euler class of the plane bundle which is orthogonal to the combing. 
To understand the extra, framing-dependent information present in the class
$c$, recall that the universal coefficient theorem for cohomology splits
$H^1(M,\Z/2)$ as
$$0 \to H^1(M,\Z) \tensor \Z/2 \to H^1(M,\Z/2) \to H^2(M,\Z) * \Z/2 \to 0,$$
and that the torsion product $H^2(M, \Z) * \Z/2$ can be viewed as the
kernel in $H^2(M,\Z)$ of multiplication by 2.
\end{proof}

Proposition~\ref{pso3action} and its proof have the following important
geometric interpretation for combings of a 3-manifold with no distinguished
framing:  Firstly, the characteristic class $c$ of a combing $b$ is defined
relative to a framing, but it only depends on the affine $H^2(M, \Z) * \Z/2$
quotient of the spin structure of the framing.  The quantity $2c$ is the
Euler class of (the normal bundle of) $b$ and is a framing-independent
cohomology class.  The maximal divisor of $c$ is framing-independent, since
it is half the maximal divisor of the Euler class. Given two combings $b_1$
and $b_2$ with characteristic classes $c_1$ and $c_2$ relative to some
framing, the difference $c_1-c_2$, the characteristic class of $b_1$ {\em
relative} to $b_2$, is also framing-independent. If the relative class is
zero and $b_1$ and $b_2$ have Hopf degrees $d_1$ and $d_2$ relative to some
framing, then $d_1-d_2$ is an element of $\Z/n$ (and not an affine
$\Z/n$-space) and is framing-independent, where $n$ is the maximal divisor
of $c_1 = c_2$.

The most important combings are those that extend to framings, namely those with
Euler class 0.  If a framing $f$ consists of a combing $b_1$ and an orthogonal
combing $b_2$, the degree of $f$ equals the Hopf degree of $b_1$, and the spin
structure comes in two pieces: a torsion product piece and a tensor product
piece.  The torsion product piece is the characteristic class of $b_1$, and
given $b_1$, the tensor product piece is determined by the twist of $b_2$ in the
normal bundle to $b_1$.  However, note that $b_1$ and $b_2$ play a symmetric
role in the definition of $f$ and that they are necessarily homotopic combings. 

It will be convenient to consider fictitious combings and framings with
fractional Hopf degree in addition to honest combings and framings. Since the
set of combings with a given characteristic class is an affine $\Z/n$-space
and the set of framings with a given spin structure is an affine $\Z$-space,
we can replace the former with an affine circle and the latter with an affine
line. A {\em fractional combing} is a formal linear combination with
non-negative real coefficients of combings of adjacent degree with total
weight 1. A {\em fractional framing} is defined similarly. In particular, the
3-sphere $\S^3$, if interpreted as $\SU(2)$, has a natural left-invariant
combing $b_L$ and a right-invariant combing $b_R$.  These differ by 1 in
degree.  The most natural combing on $\S^3$ is therefore the ambidextrous one
$$b_A = \frac12(b_L + b_R),$$
since it is invariant under an orientation-reversing diffeomorphism.
Similarly, there is a left-invariant framing $f_L$, a right-invariant
framing $f_R$, and an ambidextrous framing $f_A$.

\subsection{Combed and framed Heegaard diagrams}

Recall that a realization of a 3-manifold $M$ as a polygonal CW complex yields
a handle decomposition, where the handles are prisms that meet at rectangular
sides. The union of all rectangles is a {\em Heegaard surface}, which is a
surface that divides $M$ into two handlebodies.  In general, given a surface
$R$, a handlebody with boundary $R$ can be described by a {\em handlebody
diagram}, which is a collection of disjoint embedded circles in $R$ whose
complement is planar.  Given a handlebody diagram, the handlebody is
constructed by attaching thickened disks to one side of a thickened copy of
$R$ along the circles, and then attaching balls to the spherical boundary
components that result.  A {\em Heegaard diagram} is then a transverse pair
of handlebody diagrams on $R$, and it defines a 3-manifold for which $R$ is a
Heegaard surface.  By arbitrary convention, one of the handlebodies and its
circles is called upper and the other is called lower. Given a polygonal cell
decomposition, the perimeters of the faces are the upper circles and the
links of the edges are the lower circles of a Heegaard diagram on the
Heegaard surface of rectangles; the diagram can also be recognized as a
dissection which is dual to the tiling by rectangles:
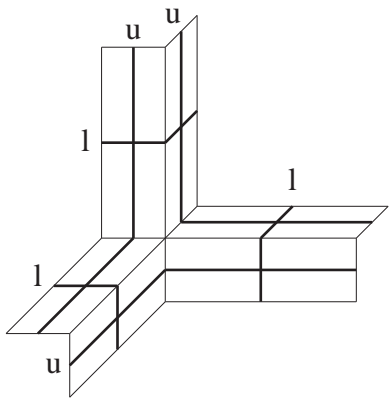
\fig
If a circle $c$ in a handlebody diagram separates two distinct components of
the complement of the diagram, then $c$ is redundant in the sense that its
removal yields another handlebody diagram which describes the same
handlebody. It follows that the only {\em minimal} handlebody diagrams on a
surface of genus $g$ are those with $g$ circles, since only those diagrams
have connected complements. Unlike in reference~\cite{Kuperberg:hopf}, most
handlebody
diagrams in this paper will be assumed to be minimal.  Note that $R$ admits
only one minimal handlebody diagram up to homeomorphism.  In this sense, only
the relative position of the two handlebody diagrams determines the topology
of a 3-manifold.

Following reference~\cite{Kuperberg:hopf}, three moves on Heegaard diagrams
suffice to render any two Heegaard diagrams for a 3-manifold $M$
equivalent.  The moves are the {\em circle slide}:
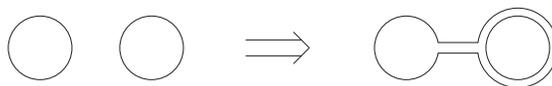
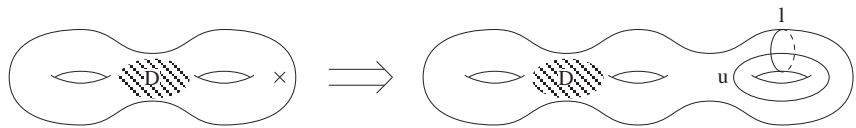
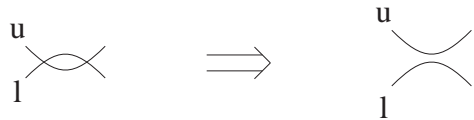
\fig
{\em stabilization}:
\fig
and {\em isotopy}:
\fig

We will need the following two well-known results concerning Heegaard
diagrams and their moves:

\begin{theorem} Any two minimal handlebody diagrams on a surface $R$ for a
handlebody $B$ are equivalent by a sequence of circle slides. \label{thslides}
\end{theorem}
\begin{proof} Let $d_1$ and $d_2$ be the two handlebody
diagrams and let $D_1$ be a collection of disjoint disks in $B$
whose boundaries form $d_1$; similarly, let $D_2$ be another such
collection whose boundaries form $d_2$.  We first show by induction
that we can eliminate the intersection of $D_1$ and $D_2$ by isotopy
and handle slides. A disk $C_1 \in D_1$ with boundary $c_1$ intersects
disks in $D_2$ in a pattern of arcs and circles:
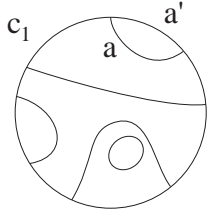
\fig
If $a$ is an innermost arc of the intersection, as shown, then the
endpoints of $a$ coincide with the endpoints of an arc $a'$ lying in
$c_1$ whose interior does not meet any circles of $d_2$.  Let $C_2 \in
D_2$ be the other disk that contains $a$ and let $c_2$ be its
boundary.  We claim that $c_2$ can be repositioned to reduce
intersection with $d_1$.  To understand how isotopy and circle slides
can be applied to $c_2$, we introduce the useful trick of cutting $R$
along every circle of $d_2$ other than $c_2$.  The result $R'$ is a
multiply-punctured torus for which $c_2$ is a meridian:
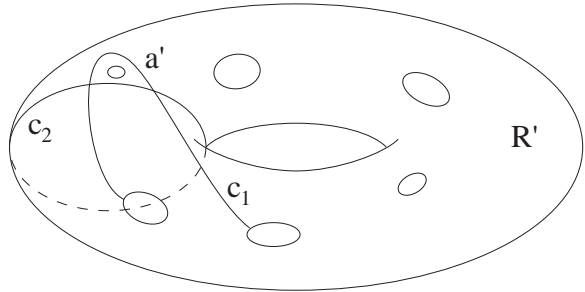
\fig
A circle slide applied to $c_2$ can be thought of as an isotopy of
$c_2$ across a puncture.  In this case, the arc $a'$ also lies in
$R'$ as indicated. (There is only one topological position for it
given that both endpoints lie on $c_2$.)  Clearly, $c_2$ can be
isotoped, perhaps across punctures of $R'$, so that an arc of it
runs parallel to $a'$:
\fig
As a result, the intersection of $c_2$ with $d_1$ decreases by at least
two points.

Thus, we can assume that $d_1$ and $d_2$ are disjoint. We show by induction
that they can be made coincident, and in this case, the hypothesis that they
describe the same handlebody is unnecessary.  Let $c_1 \in d_1$
and let $R''$ be the result of cutting $R$ along all circles of $d_2$.
The circle $c_1$ does not separate $R$, but it does separate $R''$.  Therefore
there exists a circle $c_2 \in d_2$ such that, if you glue the two
cuffs of $R''$ corresponding to $c_2$ together to obtain $R'$, $c_1$
does not separate $R'$.  The circles $c_1$ and $c_2$ are both meridians
of the punctured torus $R'$, and therefore there is an isotopy of 
$c_1$ to $c_2$, possibly involving isotopy across punctures.  Assuming
that $c_1$ and $c_2$ are coincident, we can cut $R$ along $c_1$
to obtain $R'''$, and we can apply the argument inductively
to the diagrams $d_1 - \{c_1\}$ and $d_2 - \{c_2\}$ to render
$d_1$ and $d_2$ completely parallel.
\end{proof}

\begin{theorem}[Reidemeister,Singer] Any two Heegaard surfaces for
a manifold $M$ are equivalent by a sequence of stabilization
moves. \label{thstab}
\end{theorem}
\begin{proof}(Sketch)  As we have already suggested in the context
of polygonal CW complexes, the 0- and 1-handles of
a handle decomposition of $M$ are separated from the 2- and 3-handles
by a Heegaard surface.  It is easy to realize every Heegaard surface
this way.  Following reference~\cite{Kuperberg:hopf}, it is also easy to
realize every handle decomposition by means of a Morse function on $M$.
Taking a generic path of Morse functions between any two given Morse
functions, the only transitions in the corresponding Heegaard surfaces
are stabilization moves.
\end{proof}

Although neither Theorems~\ref{thslides} and \ref{thstab} nor their
proofs refer to the isotopy move, it is implicitly allowed given the
possibility of isotopy of the lower diagram related to the upper
diagram.  The theorems show that the three given moves suffice.

If $R$ is a Heegaard surface in a 3-manifold $M$, the orientation of
$M$ induces an orientation on $R$, by the convention that a positive
tangent basis at a point for $R$ extends to a positive basis for $M$ by
appending a normal vector that points from the lower side to the upper
side.  By convention, Heegaard circles are also oriented.  Indeed, if a
Heegaard diagram comes from a polygonal complex with oriented edges and
faces, then an orientation of a face induces an orientation of its
upper circle, while an orientation of an edge induces an orientation of
its lower circle by the right-hand rule: \fig
Reversing the orientation of a circle is an obvious move that renders any two
orientations of a Heegaard diagram equivalent.

It will be important to sign-order the set of all upper and lower circles of a
minimal Heegaard diagram, where a {\em sign-ordering} of a finite set is an
orbit of the alternating group acting on the set of complete orderings of the
set. Andrew Casson pointed out the following fact to the author:

\begin{proposition} There is a canonical way to sign-order the Heegaard
circles of a minimal Heegaard diagram, if the circles are oriented.  Here
canonical means that the sign-ordering is preserved by the handle slide move
and reversed if the orientation of a Heegaard circle is reversed.
\end{proposition}
\begin{proof} Given a Heegaard surface $R$  with a minimal Heegaard diagram,
the vector space $H_1(R,\R)$ has a symplectic structure given by the
intersection form.  The upper and lower handlebodies define Lagrangian
subspaces $L_l$ and $L_u$ with respect to this form.  Viewed as homological
cycles, the circles for each handlebody form bases for $L_l$ and $L_u$, and a
sign ordering of the circles yields an orientation of the (outer) direct sum
$L_l \oplus L_u$. The vector space $H_1(R,\R)$ also has an orientation
induced by its symplectic structure.  If $L_l$ and $L_u$ are transverse,
which happens when the 3-manifold given by the diagram is a rational homology
sphere, then $L_l \oplus L_u = H_1(R,\R)$, and the canonical sign ordering of
the circles is the one that agrees with the orientation of $H_1(R,\R)$. But
this is a special case.  In the general case, the intersection form induces an
isomorphism
$$(L_l \cap L_u)^* \cong H_1(R,\R)/(L_l + L_u),$$
which means that there is an exact sequence 
$$0 \to L_l \cap L_u \to L_l \oplus L_u \to H_1(R,\R) 
\to (L_l \cap L_u)^* \to 0.$$
Recall that an orientation for all but one term of a terminating exact
sequence induces an orientation of the remaining term.  In this case, taking
either orientation of $L_l \cap L_u$ results in the same orientation of $L_l
\oplus L_u$.  This orientation yields the desired sign-ordering of the
Heegaard circles.  It is easy to check that a circle slide is a change of
basis of $L_l$ or $L_u$ with determinant 1 and that an orientation reversal
negates a basis vector, which means that the two operations respectively
preserve and reverse the sign-ordering.  \end{proof}

If we restrict the combing of a 3-manifold $M$ to a Heegaard surface $R$ and
then orthogonally project to the tangent bundle $TR$, the result is in
general a singular combing which can in principle have many types of
singularities. However, in this paper we will only consider singular combings
of $R$ with prescribed singularities in a standard position relative to a
Heegaard diagram.  We will show that such a combing canonically extends to a
combing of $M$ (or rather, it is the projection of an $M$-tangent vector
field that extends to $M$), and that these combings represent all homotopy
classes of combings. Recall that the total index of all singularities of a
combing on $R$ of genus $g$ is the Euler characteristic $2-2g$, that a
singularity of index $-1$ has the geometry
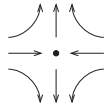
\fig
and that a singularity of index $+2$ has the geometry
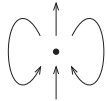
\fig
Define a {\em combing} of a minimal Heegaard diagram on $R$ to be a combing of
$R$ with $2g$ singularities of index $-1$, one on each circle, and one
singularity of index $+2$ disjoint from all circles.  The singularity of
index $-1$ on a given circle, which is called the base point of the circle,
should not lie on a crossing and the two outward-pointing vectors should be
tangent to the circle.

The extension of a combing of a diagram for $M$ to a combing of $M$ will
rely on a simple principle of interpolation to avoid discontinuities:
Consider $X \times [0,1]$ embedded in $M$ for some 1- or 2-manifold $X$.
Suppose that the $[0,1]$ fibers map to very short line segments, so
that if $p \in X$, we may think of the tangent space to $M$ at
$(p,t)$ as being the same vector space for different $t$.  Let $v$
be a field of unit vectors tangent to $M$ defined on $X \times \{0,1\}$,
and suppose further that $v(p,0)$ is never antiparallel to $v(p,1)$.
Then $v$ may be continuously extended to $X \times [0,1]$ by setting
$v(p,t)$ to be a geodesic path from $v(p,0)$ to $v(p,1)$ on the
unit sphere of possible values of $v(p,\cdot)$.

The extension of $b$ to $M$ proceeds in five steps:  Firstly, at a
singularity of index $-1$, necessarily a base point of a circle,
desingularize $b$ by making it point towards the disk attached to the circle.
Secondly, on each disk which attaches to an upper or lower circle $c$,
let $b$ be a fan emanating from the base point: \fig
If this fan were extended all the way to $c$, $b$ would be discontinuous,
because except in a neighborhood of the base point, $b$ is tangent
to the Heegaard surface $R$, while the fan is normal.  To maintain continuity,
interpolate as above between the fan and $b$ as defined on $R$.
At the point, the region on which $b$ is undefined consists of
the sole upper ball $B_u$ which is attached to $R$ and
the upper disks, the sole lower ball $B_l$ attached to $R$ and the 
lower disks, and a region which connects them
around the singularity of index $+2$, which has not been desingularized.
Geometrically, $B_u$ and $B_l$ are manifolds with corners,
because they are handlebodies with sawed handles: \fig
The third step consists of smoothing these round corners and extending
$b$ correspondingly; note
that the previous interpolation between the fan on each disk
and $b$ on each circle does not produce a situation in which
$b$ is an inward-pointing normal to the smoothed balls $B'_u$ and $B'_l$: \fig

The fourth and fifth steps are a little more difficult to visualize (which is
not to say that the first three steps are easy in this respect).  Fourthly,
identify a small topological ball $B_2$ around the singularity of index $+2$
which is geometrically a torus union a cylinder that plugs the hole of the
torus, or the surface of revolution of a barbell shape in the plane: \fig
Let $b$ be the outward normal on both the upper and lower flat disks
of $B_2$, and extend $b$ as a pair of 3-dimensional fans on
the balls $B'_u$ and $B'_l$ outside of $B_2$: \fig
As in the second step, extending $b$ all the way to the boundary
would result in discontinuities, but we may again interpolate
as we have checked that $b$ is never an inward normal on $\bd B'_u$
or $\bd B'_l$, while the fans are outward normals near the boundaries.

Fifthly, we extend $b$ to $B_2$.  For this purpose, it is best to
invert $B_2$ to recognize it as the complement of a cylinder $B'_2$ 
in the 3-sphere $\S^3$: \fig
If we coordinatize $\S^3 = \R^3 \cup \{\infty\}$ by Cartesian coordinates
$x$, $y$, and $z$ and let $B'_2$ be the cylinder given by $x^2 + y^2 \le 1$
and $|z| \le 1$, then the restriction of $b$ to $\bd B_2 = \bd B'_2$ can be
given by a simple formula, for example $(b_x,b_y,b_z) = (1,0,-z)$.
Indeed, $b$ on $\bd B'_2$ extends to a vector field $b'$ on $B'_2$ by
the same formula. We extend $b'$ to $B_2 \subset \S^3$ so that
$b'$ is the ambidextrous combing of $\S^3$; we then extend
$b$ to $B_2 \subset M$ so that it agrees with $b'$ on $B_2$.
This completes the extension of $b$ to $M$.

Given a combing $b_1$ of $M$, an orthogonal combing $b_2$ is a section of the
circle bundle orthogonal to $b_1$, and a section a circle bundle defined on a
2-skeleton of a cell complex always extends uniquely up to homotopy to the
entire cell complex.  Therefore, to describe a framing $(b_1,b_2)$
combinatorially, it suffices to describe $b_1$ as a diagram combing and then
to describe $b_2$ on the Heegaard surface $R$ and on all upper and lower
disks.  Unlike $b_1$, $b_2$ need not be in any special position in its
combinatorial description; it need only be in general position relative to
geometric objects such as crossings and upper and lower circles.  Twist
fronts indicate the position of $b_2$, where a {\em twist front} is an
arc along which $b_2$ is normal to $R$ and points from lower to upper
handlebody.  A twist front is transversely oriented in the direction
that $b_2$ rotates by the right-hand rule relative to $b_1$, and
transverse orientation is indicated by the symbol for cold
fronts on weather maps:
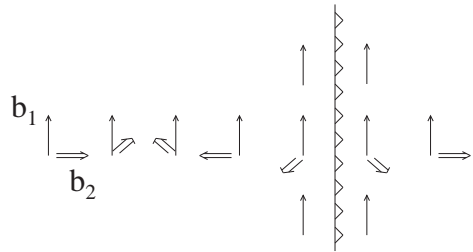
\fig
Twist fronts terminate at base points.  If the viewer's eye is
in the upper handlebody, a twist front points counterclockwise
around an upper base point and clockwise around a lower base point:
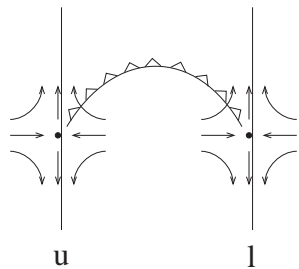
\fig
By contrast, a twist front cannot end at the singularity of index
$+2$ by conservation of twist front ends.  In essence, a $+2$
singularity is a pair of $+1$ singularities of opposite type with
a twist front connecting them, similar to the above figure.

Given a combing $b_1$ and given an orthogonal combing $b_2$ defined on $R$,
$b_2$ may or may not extend orthogonally to upper and lower disks. Since the
disks are 2-dimensional, if $b_2$ does extend, it extends uniquely up to
homotopy.  For each upper or lower circle $c$, define $\theta(c)$ to be the
total counter-clockwise rotation, in units of $1 = 360^\circ$, of the tangent
to $c$ relative to $b_1$ going around $c$ in the direction of its
orientation.  Let $\phi(c)$ be the total right-handed rotation of $b_2$
about $b_1$ going around $c$ in the direction of its rotation.
In units of $1 = 360^\circ$, $\phi(c)$ is naturally a half-integer
because of a fractional contribution at the base point.  The number $\phi(c)$
can be computed as the the number of twist fronts that cross
$c$ positively minus the number that cross negatively, with the
base point counting half as much.  An analysis of the
extension of $b_1$ to the Heegaard disk attached at $c$ shows that $b_2$
extends if and only if $\phi(c) = -\theta(c)$ when $c$ is an upper circle and
$\phi(c) = \theta(c)$ when $c$ is a lower circle.

To complete the combinatorial definition of combings and framings, we consider
a set of elementary moves on them.  All moves on uncombed minimal diagrams
are allowed, except that a circle cannot be isotoped across the $+2$
singularity:  \fig
In particular, a $+2$ singularity inside an eyelet: 
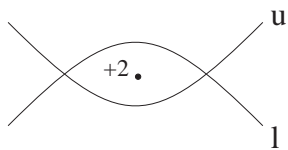
\fig
prevents the usual isotopy move.  This is restricted isotopy of
Heegaard circles.  The other moves, stabilization and a circle slides, are
unrestricted.  In addition, there are two new moves on the combings and
framings that do not affect the underlying Heegaard diagram:
\begin{itemize}
\item A {\em base point isotopy} move: \fig
The diagram indicates an isotopy of a base point on a lower circle past an
upper circle, but we also consider the isotopy of base points of upper
circles.  Also, in this diagram and later ones, a point labelled by either -1
or 2 indicates a singularity of the corresponding index.
\item A {\em base point spiral} move: \fig
\end{itemize}
Note also that homotopy of a diagram combing or framing is an allowed
operation which must be respected when defining any combed or framed
invariant, even though it is not a combinatorial move.  Homotopy
of the second combing of a framing can have a combinatorial effect
on twist fronts, in the form of four possible moves:
\begin{itemize}
\item {\em Two-point isotopy}: \fig
\item {\em Three-point isotopy}: \fig
\item {\em Circle birth}: \fig
\item The {\em exchange move}: \fig
\end{itemize}
Since only homotopy-invariant information will be used in the definition of
invariants, these moves will automatically be respected.  In particular, we
will not explicitly use the four combinatorial moves on twist fronts, but it
is a useful for practical computations to know that they suffice to reproduce
any homotopy of a second combing. Finally, a geometric analysis shows that a
spiral move does not preserve a combing (or a framing), but rather it changes
its degree by 1. More precisely, we can take as a convention that a clockwise
spiral move on an upper circle decrements the degree.  In this case a
clockwise spiral move on a lower circle increments the degree.

Before showing that these elementary moves suffice, we first show that all
combings and framings of a 3-manifold can be realized combinatorially.

\begin{lemma} Given a minimal Heegaard diagram $D$ of a 3-manifold $M$ of
genus at least 1, every combing of $M$ is realized by some combing of $D$,
and the same is true of framings. \label{lcrealize}
\end{lemma}
\begin{proof}
It is easy to find at least one combing of $D$; let $b$ denote such a combing
as well as its extension to $M$.   Recall if $b'$ is another combing, $b$ and
$b'$ have a relative characteristic class in $H^2(M,\Z)$, and if that class
is zero, they have a relative Hopf degree in $H^3(M,\Z) \cong \Z$. 
Geometrically, the relative characteristic class is the Poincar\'e dual of
the antipodal intersection of $b'$ with $b$.  \Ie, let $-b$ be $b$ with its
vectors negated, and consider $-b$ and $b'$ as sections of the unit tangent
bundle of $M$ in general position.  Then the intersection $-b \cap b'$,
projected to $M$ from $TM$, is a collection of oriented curves which
represents a homology class in $H_1(M,\Z)$ which is independent of homotopy
of $b$ and $b'$.  Its  Poincar\'e dual is the relative characteristic class.

Suppose that $b'$ is a combing of $D$ whose singularities are coincident with
those of $b$.  Since the angle between $b'$ and $b$ extends continuously to the
singularities, it defines a map $F:R \to \S^1$ from the Heegaard surface to the
circle.  If $b$ and $b'$ are in general position, their anti-parallel
intersection in $R$ yields an element $c$ in $H_1(R,\Z)$ which is Poincar\'e
dual to the cohomology class $F^*([\S^1])$.  Since there are no restrictions on
the map $F$, every element of $H_1(R,\Z)$ is realized.  Moreover, the extension
of $b$ and $b'$ to $M$ excludes the possibility that they are anti-parallel
anywhere other than on $R$. Therefore their anti-parallel intersection in $M$ is
$i_*(c)$, where the map $i_*:H_1(R,\Z) \to H_1(M,\Z)$ is induced by inclusion. 
Since this map is a surjection, every element in $H_1(M,\Z) \cong H^2(M,\Z)$ is
realized as a relative characteristic class as $b'$ varies.

The Hopf degree of $b'$ can also be set arbitrarily, because
the spiral move changes it by 1.  In conclusion, all combings of
$M$ are realized by combings of $D$.

Since $b'$ can represent any combing of $M$, then in particular, it can
represent the first combing $b_1$ of a framing $f$.  At no point in the
combinatorial definition of the second combing $b_2$ is there any restriction
of its homotopy class, and therefore all framings are realized
combinatorially also.
\end{proof}

It is useful to explicitly identify the homotopy
class of $b_2$ if it is defined combinatorially.  Given two combings $b_2$
and $b_2'$ orthogonal to $b_1$, the homological difference of their twist
fronts defines a homology class $c$ in $H_1(R,\Z) \cong H^1(R,\Z)$.  The
combinatorial constraints on twist fronts imply that $c$ lies in the image of
the injection $i^*:H^1(M,\Z) \to H^1(R,\Z)$.  (Alternatively, this can be
seen directly by defining the angle between $b_2$ and $b_2'$ and noting that
it extends to a function from $M$ to $\S^1$.) In turn, as we have already
seen, the spin stucture of $f = (b_1,b_2)$ differs from that of $(b_1,b'_2)$ by 
$$(i^*)^{-1}(c) \bmod 2 \in  H^1(M,\Z) \tensor \Z/2 \subseteq H^1(M,\Z/2),$$
while the Hopf degrees are equal.

Finally, we show that the given moves on combings and framings
render two different realizations of the same combing or framing equivalent.

\begin{lemma} Given a surface $R$ with a marked point $p$, any two minimal
handlebody diagrams for a handlebody with boundary $R$ are equivalent by isotopy
that avoids $p$ and handle slides.
\end{lemma}
\begin{proof}
The lemma is a corollary of Theorem~\ref{thslides}; it does not replace it.
Let $d$ be a minimal diagram on $R$ and let $C$ be a circle of $d$.  We
reproduce an isotopy of $C$ across the marked point $p$ by allowed moves. 
Let $R'$ be the result of cutting $R$ along all circles of $d$ other than
$c$.  Since $d$ is minimal, $R'$ is necessarily a multiply punctured torus
with $c$ a non-separating circle in $R'$.  Since the cuffs of the punctures
of $R'$ are other circles of $d$, an isotopy of $C$ across such a puncture is
effected by a circle slide.  Therefore, if $C$ is on one side of $p$, it can be
moved around $R'$ to the other side of $p$ by means of allowed operations.
\end{proof}

\begin{lemma} Any two combings of a Heegaard diagram $D$ which extend
to the same 3-manifold combing up to Hopf degree are equivalent by
base point isotopy moves and base point spiral moves.  Any two
diagram framings which yield the same 3-manifold framing up to
degree are equivalent by the same moves. \label{lcequiv}
\end{lemma}
\begin{proof}
Suppose that $b_1$ and $b_2$ are two combings of $D$ with vanishing relative
characteristic class in $M$.  Using base point isotopies, move the base
points of $b_2$ so that they coincide with those of $b_1$.  Following the
analysis in Lemma~\ref{lcrealize}, the anti-parallel intersection
of $b_1$ and $b_2$ is a class $c \in H_1(R,\Z)$.  The dual of the relative
characteristic class is $i_*(c) = 0$, where $i:R \to M$ is inclusion.
The kernel of $i_*$ is generated by the Heegaard circles of $D$, interpreted
as homology cycles on $R$.  At the same time, an isotopy of a base point
of $b_2$ all the way around a circle $d$ in $D$ changes the class $c$
by the homology class represented by $d$.  Therefore some
collection of base point isotopies renders the intersection $c = 0$.
Given $c = 0$, there is a homotopy between $b_1$ and $b_2$ that does
not move the base points; however, this homotopy may involve base point
spiral moves.

If $f_1$ and $f_2$ are two framings of $D$ which induce the same
spin structure on $M$, then firstly their first combings have vanishing
relative characteristic class. Therefore, by the previous argument, they
may be arranged so that their first combings coincide.  After this
operation, the ratio of the second combings of $f_1$ and $f_2$ is
a map $R \to \S^1$ that represents a cohomology class $c' \in H^1(S,\Z)$.
The constraints on second combings mean that $c$ lies in the image
of $i^*:H^1(M,\Z) \to H^1(R,\Z)$. Furthermore, since 
$f_1$ and $f_2$ have the same spin structure, the class
$$(i^*)^{-1}(c') \bmod 2 \in  H^1(M,\Z) \tensor \Z/2 \subseteq H^1(M,\Z/2)$$
vanishes.  In other words, $c' \in 2i^*(H^1(M,\Z))$.

The problem is to apply base point isotopies to $f_2$ to make $c'$ vanish. The
strategy for doing so is typefied in the example of $M = \S^1 \times S^2$
presented with its genus 1 Heegaard surface.  In this case, we can say by
convention that a base point isotopy all the way around the upper circle
increases both $c$ (defined using the first combings) and $c'$ by 1.  A base
point isotopy around a lower circle in the opposite direction decreases $c$
but increases $c'$. The two operations together do not change the first
combing but change $c'$ by 2, or any even number if the operations are
repeated, as desired. In general, the upper and lower circles generate
Lagrangian subgroups $L_u$ and $L_l$ of $H_1(R,\Z)$. If we identify
$H_1(R,\Z) \cong H^1(R,\Z)$ by Poincare duality, we can say that $L_u$ and
$L_l$ are subgroups of $H^1(R,\Z)$ also, with $L_u \cap L_l =
i^*(H^1(M,\Z))$.  A general collection of base point isotopies around circles
is represented by a pair $(c_u,c_l) \in L_u \oplus L_l$. These isotopies
change the homology class $c$ by $c_u+c_l$ but $c'$ by $c_u-c_l$.  Therefore
taking $c_u = -c_l \in L_u \cap L_l$, we can find base point isotopies to
make $c'$ vanish without changing $c = 0$.
\end{proof}

\section{Hopf algebras}

This section is a review of a number of results about finite-dimensional Hopf
algebras \cite{Larson-Radford,Radford:antipode,Radford:trace}.  Although the
results are not new, the arguments here separate axiomatic manipulation from
concrete considerations about finite-dimensional associative algebra, thereby
clarifying and in some cases slightly extending the results.

It will be convenient to use {\em arrow notation} for computations with
tensors \cite{Kuperberg:hopf}, which is just a graphical version of index
notation.  Given a finite-dimensional vector space $V$ and a tensor $T \in
V_1 \tensor V_2 \tensor \ldots \tensor V_n$, where each $V_i$ is either $V$
or $V^*$, we denote $T$ by its letter together with an incoming arrow for
each tensor factor of $V$ and an outgoing arrow for each tensor factor of
$V^*$.  Each arrow represents a specific tensor factor, and the permutation
of tensor factors is denoted by permuting their free ends.  For example, the
equation
\arrowpic
means that $g$ is a symmetric bilinear form on $V$.  Just as with
index notation, the tensor product of two tensors is denoted by
juxtaposing them; for example, \arrowpic
denotes $a \tensor a \in V \tensor V$, where $a \in V$.  A contraction
of two tensor factors $V$ and $V^*$ is denoted by joining arrows
head to tail; for example, \arrowpic
denotes the vector $L(v) \in V$ for a linear transformation $L$ and
a vector $v$.  Finally, the two diagrams \arrowpic
denote the identity linear transformation and $\dim V$ (the trace
of the identity), respectively.


Arrow notation is equally valid in any other pivotal, symmetric tensor
category, by which we mean a category in which tensor products and duals have
all of the usual properties that they do in the category of
finite-dimensional vector spaces. More specifically, a tensor category is a
category with an associative tensor product operation $\tensor$ on objects
and morphisms; the category is {\em symmetric} if there is a canonical
isomorphism $V \tensor W \cong W \tensor V$ that yields an action of the
symmetric group on $V^{\tensor n}$; and it is {\em pivotal} if there is a
canonical isomorphism $V^{**} \cong V$ \cite{Barrett-Westbury:spherical}. For
example, $V$ might be a (finite-dimensional) $\Z/2$-graded vector space or
{\em super}-vector space.  In this instructive case, contractions and
permutations of tensor factors obey sign rules.  For example, $\dim V$, which
is the value of an oriented circle, must be taken as the graded dimension,
the dimension of the even part minus the dimension of the odd part.  Since
tensors in arrow notation can also be interpreted as morphisms in the
category, they must be even-graded for the notation to make sense. If the
notation were extended to odd-graded tensors, the sign of the value of a
diagram would depend on the permutation sign of an ordering of its odd-graded
letters.  This sign ambiguity results from the difference between the
``internal Hom" between two graded vector spaces, which consists of all
linear transformations, and the ``external Hom", which consists of
even-graded linear transformations only.

Arrow notation partially extends to tensor categories which are not pivotal,
\eg, the category of finite- and infinite-dimensional vector spaces.  In such
categories a diagram in arrow notation is still well-defined if it is {\em
acyclic}, \ie, if there is no closed loop in the diagram in which all arrows
point in the same direction along the loop.  In the discussion below,
derivations involving acyclic diagrams are valid in these categories also.

Indeed, the invariant $\#(M,H)$ generalizes to the case in which $H$ is a Hopf
object in a tensor category in which addition is not defined.  For example,
reference~\cite{Kuperberg:hopf} shows that if $H$ is a universal, involutory
Hopf object in a universal tensor category, then $\#(M,H)$ is a complete
invariant of closed, oriented 3-manifolds.  On the other hand, we do not
consider braided tensor categories, which are categories in which the
diagrams are embedded in three dimensions and the arrows may be knotted or
linked.  Braided categories appear in the definition of the
Reshetikhin-Turaev link and 3-manifold invariants and they are implicit in
$\#(M,H)$ via the quantum double, but they are not used in any direct way in
the definition of $\#(M,H)$.  Henceforth we will use the phrase ``tensor
category'' to mean a pivotal, symmetric tensor category unless explicitly
stated otherwise.

A {\em Hopf object} $H$ (in particular a finite-dimensional Hopf algebra)
consists of two tensors \arrowpic
called multiplication and comultiplication, respectively, that
satisfy several axioms.  The axioms postulate the existence of three other
tensors: the unit $i$, the co-unit $\epsilon$, and the antipode $S$.
Multiplication is associative and unital: \arrowpic
as is comultiplication: \arrowpic
The two tensors are related by the bialgebra axiom: \arrowpic
and the axiom of the antipode: \arrowpic

Group algebras were already mentioned as examples of Hopf algebras
in Section~\ref{s:invol}.  In a group algebra, 
$$\Delta(g) = g \tensor g$$
and $S(g) = g^{-1}$ for $g$ a group element.  An important example of a Hopf
object in the graded category is an exterior algebra.  The exterior algebra
$\Lambda^*(V)$ over a finite-dimensional vector space $V$ is $\Z/2$-graded by
degree and its multiplication structure is given by the wedge product. 
Comultiplication is generated by the relation
$$\Delta(v) = v \tensor 1 + 1 \tensor v$$
for all $v \in V$.

The following lemma, proved in reference~\cite[Lemma~3.2]{Kuperberg:hopf},
establishes some basic properties of Hopf objects:

\begin{lemma} The following identities hold in any Hopf object:
\begin{itemize}
\item The tensors \arrowpic
called {\em ladders}, are inverses.  \Eg, \arrowpic
\item The counit is a multiplication homomorphism and the unit
is a co-multiplication homomorphim: \arrowpic
\item The antipode is a multiplication and co-multiplication
anti-endomorphism: \arrowpic
\item The antipode fixes the unit and the co-unit: \arrowpic
\label{lbasic}
\end{itemize}
\end{lemma}

Lemma~\ref{lbasic} uses the following shorthand for reversed multiplication
and co-multiplication: \arrowpic
In addition, the following abbrevations for multiplication or comultiplication
of more than two things will be useful: \arrowpic

\begin{lemma} The contraction $\epsilon(i)$ is 1: \arrowpic
\end{lemma}
Before proving the lemma, we explain its meaning.  The proof actually
shows that \arrowpic
For vector spaces over a field, or in any reasonable category, it is possible
to contract the identity tensor with some vector and some dual vector to
obtain 1, which would establish the lemma as stated.  Indeed, taking the
tensor product with $\epsilon(i)$ has no effect on any diagram which has an
arrow in it.  But in the abstract setting of a tensor category, it is
possible that there are scalars which cannot be expressed as contractions,
and then $\epsilon(i)$ is not necessarily 1 in the commutative semigroup of
scalars.  However, we can pass to an equally useful subcategory consisting of
the ideal generated by $\epsilon(i)$; since $\epsilon(i)$ is an idempotent,
it will be 1 in the subcategory.  This subcategory must contain
all diagrams with at least one arrow; in particular, it
contains $M$, $\Delta$, etc.

\begin{proof}
We compute: \arrowpic
\end{proof}

{\em Left integrals} $\mu_L$, {\em right integrals} $\mu_R$, {\em left
co-integrals} $e_L$, and {\em right co-integrals} $e_R$ are tensors that satisfy
the equations
\arrowpic
Traditionally, integrals are dual vectors and co-integrals are vectors,
but it will be convenient to allow them to be tensors of arbitrary type.
For example, a tensor $T$ that satisfies \arrowpic
is also a left integral.

\begin{lemma} (Existence of integrals) The tensor \arrowpic
is both a right integral and a right co-integral and has trace 1.
\label{leintegral}
\end{lemma}
\begin{proof}
Applying a ladder to \arrowpic
we obtain \arrowpic
The fact that $P_R$ is a right integral follows by applying
the inverse ladder given by Lemma~\ref{lbasic}.  The argument
that it is a right co-integral is the same with $M$ and $\Delta$
switched and with arrows reversed.  The final claim that $\Tr(P_R) = 1$
is left as an exercise.
\end{proof}

\begin{lemma} Given a right integral $\mu_R$ and a right co-integral
$e_R$, \arrowpic
\label{lsexpression}
\end{lemma}
\begin{proof}
Applying a ladder, we compute: \arrowpic
\end{proof}

\begin{corollary} (Uniqueness of integrals) Given a right integral $\mu_R$
and a right co-integral $e_R$, \arrowpic \label{cuintegral}
\end{corollary}
The corollary follows from Lemma~\ref{lsexpression} by replacing the ladder
on the left side by its inverse on the right side.

Lemma~\ref{cuintegral} suggests that $P_R$ has a rank 1, \ie, that it factors as
$e_R \tensor \mu_R$ for some $e_R$ and some $\mu_R$.  In the category of vector
spaces, this is indeed the case, but in a general tensor category, there is a
more subtle conclusion.  Applying Corollary~\ref{cuintegral} twice, we obtain
\arrowpic
(The first equality is given by taking $e_R = \mu_R = P_R$ in
Corollary~\ref{cuintegral}; the second is given by taking $e_R \tensor \mu_R
= P_R^2$.)
Let $\sigma = \Tr(P_R^2)$ be the scalar factor that appears. This is the
transposition relation.  Applying the transposition relation twice and using
$\Tr(P_R) = 1$, we obtain $\sigma^2 = 1$.  If it were the case that
$\sigma = 1$, the transposition relation would say that the inward and
outward arrows of copies of $P_R$ can be permuted separately, and it would
therefore be valid to substitute
\arrowpic
for
\arrowpic
in any diagram in which $P_R$ appears.  (If necessary, we enlarge the tensor
category in which we are working to factor $P_R$.)  However, $\sigma = -1$ in
many Hopf super-algebras, for example $\Lambda^*(V)$ when $V$ is
odd-dimensional.  When $\sigma = -1$, the transposition relation implies that
$P_R = e_R \tensor \mu_R$ for some vectors $e_R$ and $\mu_R$, but it also
implies that these vectors are odd-graded and are therefore not valid
(external) morphisms in the category of odd-graded vector spaces.
Nevertheless, it will be useful to extend arrow notation to factor $P_R$ in
the general case.  We define the notion of a formally odd-graded tensor:  A
diagram with a collection of odd-graded tensors is well-defined provided that
the tensors are sign-ordered; if the sign ordering is switched, a factor of
$\sigma$ arises. The tensor $P_R$ can be substituted by odd-graded tensors
$e_R$ and $\mu_R$, and if there are the same number of each tensor, the
sign-ordering can be indicated by dashed lines from each factor of $e_R$ to
each factor of $\mu_R$, meaning that they are ordered in such a way that the
factor of $\mu_R$ immediately follows the factor of $e_R$ that matches it,
and that
\arrowpic
As a rule, the sign-ordering can be inferred from context, and therefore it
will sometimes be omitted below.  By factoring $P_R$ and using the fact that
its trace is 1, Lemma~\ref{leintegral} can be rephrased as saying that there
exists a right integral $\mu_R$ and a right co-integral $e_R$ such that
\arrowpic
while Corollary~\ref{cuintegral} can be rephrased as
\arrowpic
for any other right integral $\mu'_R$ and right co-integral $e'_R$.

Lemma~\ref{leintegral} and Corollary~\ref{cuintegral} generalize to
left integrals and left co-integrals by symmetry, taking \arrowpic
as a definition of $P_L$.  {\em A priori}, this raises the possibility
of two sign elements $\sigma_L$ and $\sigma_R$ and separate gradings
for left and right integrals.  Happily, $\Tr(P_L^2) = \Tr(P_R^2)$ by
the following computation: \arrowpic
Therefore either trace can be called $\sigma$ and only one sign ordering
is needed for all integrals and cointegrals.  Unlike in the right integral
case, we do not define $\mu_L$ and $e_L$ so that
$P_L = \mu_L \tensor e_L$.  Rather, the two quantities are proportional
(as they must be by Corollary~\ref{cuintegral}) in a way
that will be convenient for normalization purposes.

Recall that a {\em group-like} element in a Hopf object is a vector $g$
such that \arrowpic
By Lemma~\ref{lbasic}b, $S(g)$ is also group-like: \arrowpic
It follows from the axiom of the antipode that the antipode $S$ is an
involution on group-like elements \arrowpic
In a mixture of arrow notation with more standard algebraic notation,
we will use $g^{-1}$ to mean $S(g)$ and $gh$ and $g^n$ to mean \arrowpic
for group-like elements $g$ and $h$ and $n$ an integer.

The tensor \arrowpic
is a right integral, because \arrowpic
By Corollary~\ref{cuintegral},\arrowpic
Moreover, \arrowpic
This means that \arrowpic
is a group-like element.  The vector $a$ is the {\em phase element}, and it
generates the {\em modular subgroup} of the Hopf object $H$.  We similarly
define the dual phase element $\alpha$ as \arrowpic
It is group-like in $H^*$.  We define \arrowpic
By Lemma~\ref{lsexpression}, $\mu_L$ is a left integral: \arrowpic
and since $a^{-1}$ exists, multiplication by $a$ is invertible.  The vector
$e_L$ is similarly a left co-integral.  We define
\arrowpic
In general, $\alpha^n(a^k) = q^{nk}$.  We verify that
\arrowpic
(Compare with the relation $\mu_R(e_R) = 1$.)  Also,
\arrowpic
and
\arrowpic
In conclusion,

\begin{lemma}  The vector $a$ and the dual vector $\alpha$
are group-like, the vector $e_L$ and the dual vector $\mu_L$
are a left integral and a left co-integral, and $\mu_L(e_L) = q^{-1}$.
\end{lemma}

We turn to the subject of Hopf object dualities.  The most
important such duality, which we have already
used, is the one that switches
$H = (M,\Delta)$ with its dual $H^* = (\Delta,M)$.  It is convenient
to realize this duality by reversing arrows and reflecting diagrams
about an axis, usually vertical:
\arrowpic
In addition, it is easy to check that the pairs $(M^{\op},\Delta)$,
$(M,\Delta^{\op})$, and $(M^{\op},\Delta^{\op})$ form unital, counital
bialgebras, and that the third pair is a Hopf algebra with antipode
$S$.  Call these objects $H^{\op}$, $H^{\cop}$, and $H^{\op,\cop}$,
respectively.  If either of $H^{\op}$ or $H^{\cop}$ were a Hopf algebra
with antipode $S'$, then $S' = S^{-1}$, by the following derivation:
\arrowpic

If we set \arrowpic
then $S'$ is the antipode for $H^{\op}$ and $H^{\cop}$ by the same
argument as that of Lemma~\ref{lsexpression}. Therefore,

\begin{lemma} The antipode $S$ is invertible.
\end{lemma}

Another formula for $S'$ is \arrowpic
and $S$ has a second formula as well:
\arrowpic

Note that the existence of integrals and co-integrals is the first result
which depends on the fact that the tensor category is pivotal.  
Contrariwise, suppose that $H$ is defined in a non-pivotal category but that
it has a left integral $\mu_L$ and a right co-integral $e_R$ such that
$\mu_L(e_R) = 1$.  Then since:
\arrowpic
any arrow in a diagram can be replaced by a tensor in which there is no
directed path from the tail to the head.  Therefore any diagram with cycles
can be understood as an equivalent acyclic diagram, which establishes a
pivotal structure for the tensor category containing $H$ (more precisely,
for the full subcategory whose objects are $H^{\tensor n}$.)

\begin{lemma} The antipode $S$ has the following action on 
integrals and co-integrals: \arrowpic
In particular, all integrals and co-integrals are eigenvectors of $S^2$
with eigenvalue $q$. \label{lsaction}
\end{lemma}
\begin{proof}  We compute the effect of $S$ on $e_L$:
\arrowpic
The argument is the same in the other three cases.
\end{proof}

For involutory Hopf objects with invertible dimension, the trace of
the regular representation is a right integral $\mu_R$.  In particular,
$\mu_R(AB) = \mu_R(BA)$ in this case.  In the general case, $\mu_R$ has
a property similar to that of a trace:

\begin{lemma}  The tensors $\mu_R$ and $e_R$ satisfy:
\arrowpic
\label{lqtrace}
\end{lemma}
\begin{proof} We give the proof for $\mu_R$.  Using expressions for
$S$ and $S^{-1}$ and a formula for $e_L$,
\arrowpic
Since $S$ is invertible, the tensor 
\arrowpic
is also and we can remove it from both sides.  Applying $S^{-2}$ to both
inward arrows on both sides then yields the desired result, using the fact
that $\mu_R$ is an eigenvector with eigenvalue $q^{-1}$.
\end{proof}

We define
\arrowpic
for any integer $n$.  Note that $\mu_{1/2} = \mu_L$, $\mu_{-1/2} = \mu_R$,
$e_{1/2} = e_L$, and $e_{-1/2} = e_R$.

\begin{theorem}[Radford] \arrowpic
or, equivalently, $\Ad_\alpha^* \circ \Ad_a = S^4$.
\label{thradford}
\end{theorem}
\begin{proof} First express $S^2$ as \arrowpic
using both formulas for $S$ and substitute for $S^{-1}$ in the middle.
Applying $\Ad_{\alpha^{-1}}^*$, we obtain \arrowpic
Applying $\Ad_{a^{-1}}$, we obtain \arrowpic
The result follows from the fact that
the operators $\Ad_\alpha^*$ and $\Ad_a$ commute with $S^2$ (check).
\end{proof}

Theorem~\ref{thradford} has the following important corollary for
finite-dimensional Hopf algebras (as opposed to general Hopf objects): Since
distinct group-like elements are linearly independent, $a$ and $\alpha$ have
finite order.  Since $\Ad_\alpha^*$ and $\Ad_a$ commute with each other as
well as with $S^2$, it follows that $S$ has finite order.

A Hopf object $H$ is {\em balanced} if $\Ad_a = S^2$, which by
Theorem~\ref{thradford} means that $H^*$ is also balanced.  In the unbalanced
case, define
\arrowpic
The tensor $T$ is the {\em tilt} of $H$, since it measures the extent
to which $H$ fails to be balanced.

\begin{lemma}  The tilt map $T$ is an automorphism of $H$ that fixes
integrals and co-integrals.
\end{lemma}

The proof is left as an exercise.

\subsection{Quantum groups}

Let $\frak g$ be a simple Lie algebra over $\C$ and let $\frak g^+$ be a Borel
subalgebra.  The universal enveloping algebras $U({\frak g})$ and $U({\frak
g}^+)$ are Hopf algebras that admit deformations $U_q({\frak g})$ and
$U_q({\frak g}^+)$, where $q$ can be  a non-zero complex number or an
indeterminate \cite{Drinfeld,Lusztig}.   (It will appear as if $U_q({\frak
g}^+)$ depends on a choice of $q^{1/2}$, but the dependence disappears with a
slightly different choice of generators.  Moreover, our convention is
consistent with the definition $q = \alpha(a)$ given above.) Furthermore, if
$q$ is a root of unity, these Hopf algebras admit finite-dimensional
quotients $u_q({\frak g})$ and $u_q({\frak g}^+)$. Each of these four
deformations has an established role in theory of quantum topological
invariants \cite{RT:graphs,RT:manifolds}. As mentioned in the introduction,
the invariant $\#(M,u_q({\frak g}^+))$ is an important example of $\#(M,H)$. 
In light of this use, the goal of this section is the following lemma:

\begin{lemma} The Hopf algebra $u_q({\frak g}^+)$ is balanced.
\label{lfrenkel}
\end{lemma}

Before proving the lemma, we give a quick definition of $u_q({\frak g}^+)$
(which in any case varies in the literature):  Let $\frak g$ have rank $n$,
let $\alpha_1,\ldots,\alpha_n$ be the simple roots of the root system of
$\frak g$, let $(\cdot,\cdot)$ be the dual Killing form, and let $\langle
\alpha, \beta \rangle = (\alpha,\beta)/(\beta,\beta)$, so that $\langle
\alpha_i, \alpha_j \rangle$ is the Cartan matrix. Let $r$ be the order of
$q$, and choose a square root $q^{1/2}$.  Then $u_q({\frak g}^+)$ is
generated as a unital algebra by the generators $E_i$ and $K_i$ for $1 \le i
\le n$ with the relations
\begin{eqnarray*}
K_i K_j & = & K_j K_i, \\ 
K_i E_j & = & q^{(\alpha_i,\alpha_j)/2} E_j K_i, \\
\Ad_\Delta(E_i)^{1-\langle \alpha_j,\alpha_i\rangle}(E_j) & = & 0, \\ 
K_i^{2r} & = & 1, \\
E_i^r & = & 0,
\end{eqnarray*}
where $\Ad_\Delta(X)$ is the linear operator defined as:
\arrowpic
$$\Delta(E_i) = E_i \tensor 1 + K_i \tensor E_i,$$
$$S(E_i) = -K_i^{-1}E_i,$$
and
$$S(K_i) = K_i^{-1}.$$
(The trick of expressing the quantum Serre relations in terms of the
quantum adjoint was suggested to the author by Marc Rosso.)
The rest of co-multiplication is generated as an algebra homomorphism from this
relation and
$$\Delta(K_i) = K_i \tensor K_i.$$
Similarly, the rest of the antipode map is generated as an algebra
anti-automorphism.  The unit
is just written as $1$.
The counit $\epsilon$ is generated as an algebra homomorphism by
$$\epsilon(E_i) = 0,$$
$$\epsilon(K_i) = 1.$$

\begin{proof} (Sketch) The elements $K_i$ are group-like and generate a subgroup
$C$ of $u_q({\frak g}^+)$; the group algebra $\C[C]$ is a Hopf subalgebra. 
There is a sequence $i_1,\ldots,i_k$, where $k$ is the number of
positive roots of $\frak g$, such that
$$E_{\max} = E^{r-1}_{i_1} E^{r-1}_{i_2} \ldots E^{r-1}_{i_k}$$
is non-zero, but that, for any such sequence, $E_{\max}E_i = 0$ for all
$i$.  Moreover, each $i$ appears $m_i$ times as some $i_j$,
where $2\rho = \sum_i m_i \alpha_i$ is the sum of the positive
roots of $\frak g$.
Then 
$$e_L = \left(\sum_{K \in C} K\right)E_{\max}$$
is clearly a left co-integral, and similarly
$$e_R = E_{\max}\left(\sum_{K \in C} K \right)$$
is a right co-integral.  (We omit formulas for $\mu_L$ and $\mu_R$,
since we will not use them here.)  If we define $\alpha$ as usual by
$$ X e_R = \alpha(X) e_R,$$
then
$$\alpha(E_i) = 0$$
and
$$\alpha(K_i) = q^{(r-1)(\rho,\alpha_i)} = q^{-(\alpha_i,\alpha_i)},$$
since $\langle \rho,\alpha_i\rangle = 1$ for all $i$.
We now compute the effect of $\Ad_\alpha^*$, which is an algebra
automorphism, on $E_i$ and $K_i$:
$$\Ad_\alpha^*(K_i) = (\alpha \tensor I \tensor \alpha^{-1})(\Delta_3(K_i))
= \alpha(K_i) \tensor K_i \tensor \alpha^{-1}(K_i) = K_i,$$
while
\begin{eqnarray*}
\Ad_\alpha^*(E_i) & = & (\alpha \tensor I \tensor \alpha^{-1})(\Delta_3(E_i)) \\
& = & \alpha(E_i) \tensor 1 \tensor 1 + \alpha(K_i) \tensor E_i \tensor 1
+ \alpha(K_i) \tensor K_i \tensor \alpha^{-1}(E_i) = q^{-(\alpha_i,\alpha_i)} E_i.
\end{eqnarray*}
Here $\Delta_3(X)$ is the triple co-product, meaning
\arrowpic
This value agrees with
$$S^2(E_i) = K_i^{-1} E_i K_i = q^{-(\alpha_i,\alpha_i)} E_i.$$
Since $\Ad_\alpha^*$ and $S^2$ agree on generators of $u_q({\frak g}^+)$
and they are both algebra automorphisms, the two operators agree
always and $u_q({\frak g}^+)$ is balanced, as desired.
\end{proof}

\section{The invariant}
\label{sframed}

Let $M$ be a closed 3-manifold with a minimal Heegaard diagram $D$. Orient all
circles and give them the canonical sign-ordering relative to the
orientation.  Let $f = (b_1,b_2)$ be a framing of $D$.  Recall the definition
of
$\theta(c)$ and $\phi(c)$ for Heegaard circles $c$.  For each point $p$ on
some circle $c$ of $D$ with base point $o$, we also define $\theta(p)$ be the
counterclockwise rotation of the tangent to $c$ relative to $b_1$ from $o$ to
$p$, as before in units of $1 = 360^\circ$.  If $p$ is a crossing, then two
rotation angles are defined; call them $\theta_l(p)$ and $\theta_u(p)$.
Arrange the circles so that upper and lower circles are not only transverse,
but also orthogonal when they cross, so that $\theta_l(p) - \theta_u(p)$ is always
$\frac n2 + \frac14$ for some integer $n$.  Let $\phi(p)$ be the total
right-handed twist of $b_2$ around $b_1$ from $o$ to $p$, and similarly define
$\phi_u(p)$ and $\phi_l(p)$. Using twist fronts, $\phi(p)$ can be computed
as the total sign of all fronts crossed from $o$ to $p$, not counting the
front that terminates at $o$ itself.

Let $H$ be a Hopf object in a pivotal tensor category. Let $\mu_R$, $\mu_L$,
$e_R$, and $e_L$ be integrals and co-integrals of $H$ such that $\mu_R(e_R) =
\mu_R(e_L) = \mu_L(e_R) = 1$. Recall also the tensors, $a$, $\alpha$, and $T$
associated to $H$ and the scalars $\sigma$ and $q$.

Define the quantity $\#(D,H)$ as follows:  Replace each upper circle $c$
with an $M$ tensor with one inward arrow for each crossing and the outward
arrow with an integral at the base point, with the arrows ordered
as indicated:
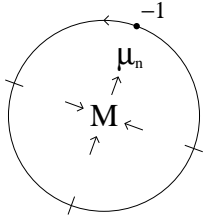
\arrowpic
Here $n = -\theta(c)$.  Replace each lower circle $c$ with a $\Delta$ tensor
with an outward arrow for each crossing and the inward arrow with a
cointegral at the base point, with the arrows ordered as indicated:
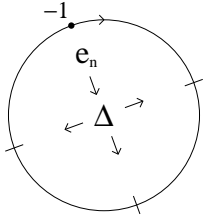
\arrowpic
Here $n = \theta(c)$.  Since there is an integral or a co-integral
for each circle, give them the same sign-ordering as that of the
circles.  Replace each crossing by the tensor:
\arrowpic
where $n = 2(\theta_l(p) - \theta_u(p)) - \frac12$, $k = \phi_l(p) -
\phi_u(p)$, and $p$ is the crossing point.  Finally, contract all tensors
corresponding to circles and crossings according to incidence.

\begin{theorem} The quantity $\#(D,H)$ depends only on $M$ and its framing.
\end{theorem}

\begin{proof}
We demonstrate invariance or covariance under each type of diagram move. For most of the
moves, we use the duality between $H$ and $H^*$ to cut the calculation in half.
\begin{itemize}
\item Orientation reversal.  If the orientation of a lower circle $c$
is reversed, $n = \theta(c)$ is subtracted from $\theta_l(p)$ for
each crossing $p$ on $c$ and then $\theta(c)$ is negated.  The crossings
on $c$ reverse order, and the sign-ordering of all circles also reverses.
Therefore a tensor for a lower circle such as:
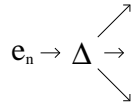
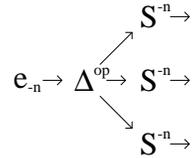
\arrowpic
is replaced with
\arrowpic
By Lemma~\ref{lsaction}, these two are equal.
\item Rotation front spiral. 
This move has the effect of changing $\phi_l(p)$ or $\phi_u(p)$ by one
for all crossings on an upper or lower circle.  A tensor for a lower circle
such as
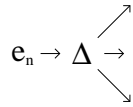
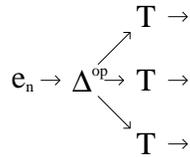
\arrowpic
might change to
\arrowpic
These are equal because $T$ is an automorphism of all of the intrinsic
structure of $H$.
\item Two-point isotopy.  After appropriate orientation reversal,
a digon which is the starting point of a two-point isotopy move
can look like:
\fig
The corresponding tensor is
\arrowpic
for some $n$, which is equivalent to the tensor after the move by the
computation shown.
\item Base point isotopy.  For simplicity, we assume a rotation of 
a lower circle $c_l$ that moves its base point past a crossing $p$ with
an upper circle $c_u$, such that the identity tensor $I = S^0T^0$ is
assigned to $p$.  Following Figure~18, the relevant
piece of $\#(M,H)$ might be:
\arrowpic
where $n = \theta(c_l)+\frac12$ and $k = \frac12-\theta(c_u)$.  The computation
shown uses Lemma~\ref{lqtrace}.  The final expression matches
the result of the rotation, because $\phi_l(p)$ changes from $0$
to $n-1$, $\theta_l(p)$ changes from $0$ to $n-\frac12$, $\theta_u(p)$
decreases by $\frac12$, and $\theta(c_u)$ and $\theta(q)$ decrease by one
for every point $q$ after $p$ on $c_u$ (and similarly $\phi(c_u)$ and
$\phi(q)$ change in the opposite direction).
\item Circle slide.  A circle slide might typically
look like this:
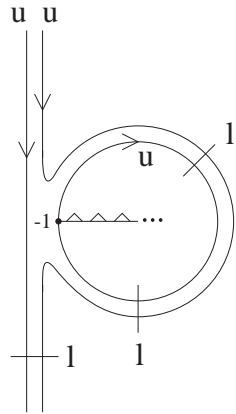
\fig
We again assume for simplicity that all crossings before the move have
tensor $I = S^0T^0$.  Then the relevant piece of $\#(M,H)$ before the move
is:
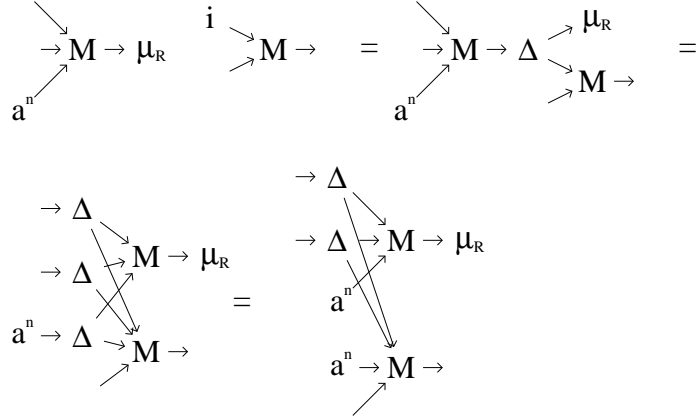
\arrowpic
where $c_1$ slides past $c_2$ and $n = \frac12 - \theta(c_2)$. The derivation
given arrives at the result of the handle slide move, except for one omitted
step.  Just as with the circle rotation move, a factor of $a^n$ appears which is
in the wrong position.  Conjugating by $a^n$ yields $S^{2n}T^{-n}$ tensors. 
These tensors arise in the move because $\theta(p)$ decreases and $\phi(p)$
increases by $n$ for all points $p$ along $c_2$ after the position of the circle
slide.  In the given example, $n = 0$.
\item The stabilization move yields a scalar factor which equals 1: \arrowpic
\item Base point spiral:  A clockwise spiral on a lower circle $c$ increases
$\theta_l(p)$ by one
at every crossing of $p \in c$.  This brings in a factor of $S^2$,
and by Lemma~\ref{lsaction}, the invariant gains a factor of $q$.
In general, changing the Hopf degree of the framing of $M$ by $n$
changes $\#(M,H)$ by $q^n$.
\end{itemize}
\end{proof}

\subsection{The combed invariant}

The definition of $\#(M,H)$ for $H$ balanced and $M$ combed is in fact
exactly the same as for $H$ arbitrary and $M$ framed. Since the tilt
map $T$ is the identity, the second combing $b_2$ of a framing $f =
(b_1,b_2)$ is irrelevant.  In this case, all constructions of
Section~\ref{sframed} are valid if $b_1$ is a combing that does not
extend to a framing, with the conclusion that $\#(M,H)$ is defined
even when the combing of $M$ has non-zero Euler class.

\section{Properties and examples}

The manifold $\R P^3$, if interpreted as the group $\SO(3)$,
has two natural combings:  The left-invariant combing and the right-invariant
combing.  These two combings differ in their characteristic class, 
but not in their Euler class; since $H^2(\R P^3,\Z)$ is entirely 2-torsion,
all combings extend to framings.  (One consequence is the well-known
fact that the left-invariant
and right-invariant framings of $\R P^3$ differ in spin structure.)
Indeed, we could just as well say that $\R P^3$ is framed rather than
combed.
A simple computation using the genus 1 Heegaard diagram of $\R P^3$
demonstrates that 
$$\#(\R P^3,H) = \Tr(S)$$
with one combing for any $H$, and that
$$\#(\R P^3,H) = \Tr(S^{-1})$$
for the other combing.  If $H = u_q({\sl(2)}^+)$, the basis
$$\{K^j (K^{-1} E^2)^k, K^j (K^{-1} E^2)^k E\}$$
is convenient, because $S$ permutes the basis elements up to scalar
factors.  The only basis elements that are (projectively) fixed
are those of the form $(K^{-1} E^2)^k$ and $K^r (K^{-1} E^2)^k$
with $0 < k < r/2$.  The result is that
$$\Tr(S) = 2{1 - q^{-\lfloor {r + 1 \over 2}\rfloor}\over 1 - q^{-1}}$$
and
$$\Tr(S^{-1}) = 2{1 - q^{\lfloor {r + 1 \over 2}\rfloor}\over 1 - q}.$$
These values are similar to values for the Reshetikhin-Turaev 3-manifold
invariants, but not quite the same.

The manifold $\S^2 \times \S^1$ has many different combings, and in
general the value of the invariant is
$$\#(\S^2 \times \S^1,H) = \mu_R(a^n)e_R(\alpha^k)$$
for some $n$ and $k$.  If $q \ne 1$, this expression must be zero.
When $n = k = 0$, $\#(\S^2 \times \S^1,H)$ reduces to $\Tr(S^2)$, which is
non-zero for Hopf algebras precisely when $H$ is semisimple and
co-semisimple.  Furthermore, in characteristic 0, it is known that
$H$ is semisimple if and only if it is involutory \cite{Larson-Radford}.

A connected sum $M_1 \# M_2$ of two 3-manifolds has a Heegaard
diagram which is also a connected sum; the Heegaard circles of the
two pieces are disjoint.  With suitable combings or framings,
$$\#(M_1 \# M_2, H) = \#(M_1,H)\#(M_2,H).$$

The symmetries in the definition of $\#(M,H)$ immediately yield
the identities:
$$\#(M,H) = \#(-M,H^{\op}) = \#(-M,H^{\cop})$$
and
$$\#(M,H^*) = \#(M,H),$$
where $-M$ is $M$ with reversed orientation.

Suppose that a combing of $M$ has divisor $d$ but $q^d \ne 1$, where
as before $q = \alpha(a)$ for some Hopf algebra $H$.   Then 
$\#(M,H) = 0$, because on the one hand, changing the degree of the
combing by $d$ changes $\#(M,H)$ by $q^d$, but on the other
hand, it does not change the combing at all.

Since $\Lambda^*(\C)$ is commutative and co-commutative in the graded sense,
$\#(M,\Lambda^*(\C))$ does not depend on the order of the crossings on each
upper and lower Heegaard circle.  Therefore $\#(M,\Lambda^*(\C))$ can only
depend on the intersection matrix between the upper and lower crossings.
Indeed, it is the determinant of this matrix up to sign; and the
sign-ordering normalization ensures that $\#(M,\Lambda^*(\C))$ is
non-negative.  Therefore
$$\#(M,\Lambda^*(\C)) = |H_1(M,\Z)|$$
when the right side is finite, and
$$\#(M,\Lambda^*(\C)) = 0$$
otherwise.

If a finite group $G$ acts by automorphisms on a Hopf algebra or
super-algebra $H$, then $G \ltimes H$ is well defined as a Hopf algebra
semi-direct product.  In particular, if $V$ is a linear representation
of $G$, then $G \ltimes \Lambda^*(V)$ is an interesting Hopf
super-algebra.  Unfortunately, $\#(M,G \ltimes \Lambda^*(V))$ is
not interesting. It expands as a sum over homomorphisms from $\pi_1(M)$
to $G$.  If a given homomorphism $\pi_1(M) \to G$ is composed with
the representation $G \to \End(V)$, the result is 
a flat vector bundle $E$ over $M$ with fiber $V$ at the base point.  The
corresponding term in $\#(M,G \ltimes \Lambda^*(V))$ is the determinant
of the middle
term of the chain complex
$$0 \longto C_3(M,E) \longto C_2(M,E)
\longto C_1(M,E) \longto C_0(M,E) \longto 0.$$
Unless the homomorphism $\pi_1(M) \to G$ is trivial, this determinant
is 0.  The determinant of the whole complex is the Reidemeister
torsion of $E$ and is an interesting invariant, but it does not
appear in $\#(M,G \ltimes \Lambda^*(V))$.
Thus, $\#(M,G \ltimes \Lambda^*(V))$ almost sums the Reidemeister
torsion of $M$ over linear representations of $\pi_1(M)$ into $V$
that factor through $G$, but instead evaluates to $\#(M,\Lambda^*(\C))$.

There are several reasons to believe that there is a
generalization of $\#(M,H)$ that involves the 0-cells and the 3-cells
of a cell decomposition of $M$ rather than just the edges and faces:
\begin{itemize}
\item There might be an analogue of the non-reduced
cohomology set $H^1(M,G)$, as well as the weighted enumeration
of $H^1(M,G)$ considered by Dijkgraaf and Witten \cite{Dijkgraaf-Witten}.
\item The invariant $\#(M,G \ltimes \Lambda^*(V))$
comes close to an interesting, classically considered invariant.
\item The invariant $\#(M,H)$ is not always a TQFT. In particular,
$\#(M,u_q({\frak g}^+))$ is not a Jones-Witten TQFT.
\end{itemize}

The treatment of $\#(M,H)$ in this paper is from the point of view
of deriving information about 3-manifolds using Hopf algebras and
Hopf objects.  An equally useful point of view, which we did
not address, is the converse:
What can we learn about Hopf algebras using 3-manifolds?
Reference~\cite{Kuperberg:hopf} establishes a noteworthy result
in this direction:  Two expressions in an abstract involutory Hopf object 
(or in the universal involutory Hopf object) are
axiomatically equal if and only if two related 3-manifolds with
boundary are homeomorphic.  It would be interesting to extend this
result to arbitrary Hopf objects.  Such an extension would involve
interpreting the axioms of a Hopf object as moves on some kind
of topological object, probably a suitably decorated 3-manifold.

One of the merits of $\#(M,H)$ in the involutory case is that
the definition is particularly simple.  Unfortunately, the
general definition of $\#(M,H)$ has been a disappointment by
comparison.  Perhaps if the axioms of a Hopf object were properly
understood topologically, it would lead to a simpler definition
of $\#(M,H)$.

\section{Acknowledgements}

We wish to especially thank Igor Frenkel for his sustained interest in this
paper and for the proof of Lemma~\ref{lfrenkel}.  We also thank Andrew Casson
and Vaughan Jones for their encouragement and early attention to the author's
work, and we thank Bruce Westbury, Thomas Mattman, and the referee for their
careful reading of various drafts of this paper.


\end{document}